\documentclass{amsart}
\usepackage{hyperref}
\hypersetup{colorlinks,citecolor=blue,linkcolor=blue,urlcolor=blue}

\usepackage{amsmath,amssymb,amsthm}
\usepackage{comment,enumerate}
\usepackage[capitalize]{cleveref}
\usepackage[all]{xy}

\crefname{equation}{equation}{equations}
\crefname{section}{Section}{Sections}
\crefname{chapter}{Chapter}{Chapters}

\newtheorem{lem}{Lemma}[section]
\newtheorem{thm}[lem]{Theorem}
\newtheorem*{thm*}{Theorem}
\newtheorem{cor}[lem]{Corollary}
\newtheorem{prop}[lem]{Proposition}

\theoremstyle{definition}

\newtheorem{question}[lem]{Question}

\newcommand{\RR}{\mathbb{R}}

\newcommand{\PP}{\mathbb{P}}

\newcommand{\cha}{\operatorname{char}}

\newcommand{\gds}{G_d^{sep}}

\newcommand{\tw}[1]{{\widetilde{#1}}}

\newcommand{\benum}{\begin{enumerate}[(1)]}
\newcommand{\eenum}{\end{enumerate}}

\numberwithin{equation}{section}
\begin{document}

\title{Analogs of the Shapiro Shapiro Conjecture in Positive Characteristic}

\author{Ryan Eberhart}
\address{Department of Mathematics, The Pennsylvania State University, University Park, PA 16802}
%\curraddr{}
\email{rde3@psu.edu}
%\thanks{}

%\subjclass[2010]{Primary 12F10; Secondary 14D99}
%\subjclass[2010]{Primary 12F10; Secondary 14D05, 14L30}

\date{}

\dedicatory{}

\begin{abstract}
Motivated by the Shapiro Shapiro conjecture, we consider the following: given a field $k$, under what conditions must a rational function with only $k$-rational ramification points be equivalent (after post-composition with a fractional linear transformation) to a rational function defined over $k$? The main results of this paper answer this question when $k$ has characteristic 2 or 3. We also show the insufficiency of several natural conditions in higher characteristic.
\end{abstract}

\maketitle

\section{Introduction}
In \cite{shapiro}, Eremenko and Gabrielov proved a special case of the Shapiro Shapiro conjecture. Their result states that a (univariate) rational function defined over the complex numbers with only real ramification points is equivalent to a rational function defined over the real numbers. Here and throughout this paper rational functions $f_1(x)$ and $f_2(x)$ are considered \textit{equivalent} if there is a fractional linear transformation $\sigma(x)$ such that $f_2=\sigma\circ f_1$. Motivated by this result, it is natural to ask the following:

\begin{question}\label{mainquest}
Fix a field $k$. Under what conditions, on for example the degree or ramification indices, must a rational function defined over $\overline{k}$ with only $k$-rational ramification points be equivalent to a rational function defined over $k$?
\end{question}

In \cite[Theorem 1.3]{fabershapiro}, Faber and Thompson provide a necessary and sufficient condition on a field $k$ of characteristic 0 so that every degree 3 rational function defined over $\overline{k}$ with only $k$-rational ramification points must be equivalent to a rational function defined over $k$. We note that the arguments in \cite{fabershapiro} also hold for fields $k$ of characteristic greater than 3, which we shall make use of in \cref{negativesect}.

Instead of restricting to low degree rational functions, we consider \cref{mainquest} over fields of positive characteristic. The primary results of this paper can be summarized as follows (the Wronskian is defined in \cref{prelimsect} below):

\begin{thm}(\cref{mainthm}, \cref{finitefieldsimple}, \cref{char2cor}, \cref{char3cor})\label{summarythm}
\begin{enumerate}[(i)]
\item Let $k$ be a field of characteristic 3. Every simply ramified (i.e.\ no ramification index is greater than two) rational function defined over $\overline{k}$ with only $k$-rational ramification points is equivalent to a rational function defined over $k$.

\item Let $k$ be a finite field of characteristic greater than 3. There exists a degree 3 simply ramified rational function defined over $\overline{k}$ with only $k$-rational ramification points which is \textit{not} equivalent to a rational function defined over $k$.

\item Let $k$ be a non-algebraically closed field of characteristic 2. Suppose that $f(x)$ is a degree $d>1$ rational function defined over $\overline{k}$ with only $k$-rational ramification points. Then there are infinitely many equivalence classes of degree $d$ rational functions with the same Wronskian as $f(x)$ which are not equivalent to a rational function defined over $k$. The same result holds if $k$ is a non-algebraically closed field of characteristic 3 when $f(x)$ is not simply ramified.
\end{enumerate}
\end{thm}
These results provide a complete answer to \cref{mainquest} in characteristic 2 and 3. In particular, \cref{summarythm} asserts that the correct condition for \cref{mainquest} in characteristic 3 is that the ramification points do not coalesce.

We briefly remark on the differences between previous results in characteristic 0 and the positive characteristic results of this paper. In \cite{shapiro} it is shown that no additional conditions are necessary in \cref{mainquest} when $k=\RR$ for a generic choice of distinct real ramification points. Using the fact that the Wronskian (defined in \cref{prelimsect} below) is a quasi-finite map, a limiting argument (see for example the discussion proceeding \cite[Theorem 1.3]{fullshapiro}) extends this result for any configuration of real ramification points. However, the Wronskian is \textit{not} quasi-finite in positive characteristic. This explains why the simply ramified result in \cref{summarythm}(i) does not imply the same result for the case when the ramification points coalesce (which would contradict \cref{summarythm}(iii)).

\section{Definitions and Preliminaries}\label{prelimsect}
Throughout this paper we shall implicitly assume that all rational functions are separable. Let $k$ be a field and $f(x)=g(x)/h(x)$ a degree $d$ rational function defined over $\overline{k}$ (i.e.\ $g(x),h(x)\in\overline{k}[x]$), where $g(x)$ and $h(x)$ are coprime polynomials. We say that $f(x)$ is \textit{simply ramified} if none of its ramification indices are greater than two. In this paper we shall use a measure of ramification related to the ramification index, called the differential length. The \textit{differential length} of $f(x)$ at a closed point $P\in\PP^1_{\overline{k}}$, denoted $l_P$, is the length of the sheaf of relative differentials as an $\mathcal{O}_{\PP^1_{\overline{k}},P}$ module. We introduce a convenient method to calculate the differential lengths below.

The \textit{Wronskian} of $f(x)$ is the polynomial $Wr(f(x))=h(x)g'(x)-g(x)h'(x)$. After factoring over $\overline{k}$ we have that
$$Wr(f(x))=\alpha\prod (x-c_i)^{l_i},$$
where each $l_i$ is positive, $\alpha\in\overline{k}$, and each $c_i\in\overline{k}$ is distinct. The $c_i$, considered as points on $\PP^1_{\overline{k}}$, are precisely the affine ramification points of $f(x)$. Each $l_i$ is equal to the differential length of $f(x)$ at $c_i$. By the Riemann Hurwitz formula the differential length at $\infty$ is equal to $2d-2-\deg(Wr(f(x)))$. Note that if $f(x)$ is tamely ramified at $P$ then $l_P$ is equal to the ramification index minus one; if $f(x)$ is wildly ramified at $P$ then $l_P$ is strictly larger than that.

Given a positive integer $d$, let $G_d^{sep}$ denote the space of separable two dimensional linear series on $\PP^1_k$. Any such linear series corresponds to an equivalence class of separable degree $d$ rational functions. Since the Wronskian of an equivalence class of rational functions is defined up to multiplication by a constant, the Wronskian induces a morphism $Wr:G_d^{sep}\rightarrow \PP^{2d-2}_k$. Here we consider $\PP^{2d-2}_k$ as the projectivization of the space of polynomials of degree at most $2d-2$. We consider $\gds$ as a scheme defined over $k$. See \cite[Section 4]{myconj} for a more detailed discussion.

Fix a positive integer $d$ and a finite tuple of positive integers $T=(l_1,...,l_n)$ such that $\sum l_i=2d-2$. Following \cite[Lemma 4.1]{myconj}, let $X_T$ denote the locally closed subscheme of $G_d^{sep}$ corresponding to equivalence classes of degree $d$ rational functions ramified at unspecified points $P_1,...,P_n$ such that $l_{P_i}$ is equal to $l_i$. As mentioned in the introduction, in positive characteristic the map $Wr:\gds\rightarrow \PP^{2d-2}_k$ is not quasi-finite. However, in certain circumstances the restriction of this map to one of the $X_T$ subschemes is quasi-finite, which will be an important ingredient in our results.

As our final preliminary, we find a unique representative for each equivalence class of rational functions which will be convenient in what follows.

\begin{lem}\label{form}
Let $k$ be a field.

\begin{enumerate}[(i)]
\item Every rational function defined over $\overline{k}$ is equivalent to a unique rational function $g(x)/h(x)$, where $g(x)$ and $h(x)$ are monic, $\deg(g(x))>\deg(h(x))$, and $g(x)$ contains no $\deg(h(x))$ term. This will be referred to as the \textit{standard form} of a rational function.

\item Such a rational function is equivalent to a rational function defined over $k$ if and only if $g(x),h(x)\in k[x]$ in its standard form. The latter condition is equivalent to the corresponding point in $G^{sep}_d$ having residue field $k$.
\end{enumerate}
\end{lem}

\begin{proof}
Let $f(x)$ be a degree $d$ rational function. We begin with (i) and first show the existence of a standard form. After post-composing $f(x)$ with a fractional linear transformation sending $f(\infty)$ to $\infty$, we have a rational function of the form $g(x)/h(x)$, where $\deg(g(x))>\deg(h(x))$. By dividing $g(x)$ and $h(x)$ by the leading coefficient of $h(x)$ we may assume that $h(x)$ is monic. Post-composing with $x\mapsto x/a$, where $a$ is the leading coefficient of $g(x)$, yields a form where $g(x)$ is also monic. Finally, post-composing with $x\mapsto x-b$, where $b$ is the coefficient of $x^{\deg(h(x))}$ in $g(x)$, leads to the desired form.

To show uniqueness, it suffices to prove that distinct standard forms cannot be equivalent. To this end, suppose that $f_1(x)=g_1(x)/h_1(x)$ and $f_2(x)=g_2(x)/h_2(x)$ are two standard forms such that $f_2=\sigma\circ f_1$, where $\sigma(x)$ is a fractional linear transformation. Since $\deg(g_i(x))>\deg(h_i(x))$ and each $g_i(x)$ and $h_i(x)$ is monic, it must be the case that $\sigma(x)=x+a$, where $a\in\overline{k}$. Since $f_1(x)$ and $f_2(x)$ are the same degree and each fixes $\infty$, we have that $\deg(g_1(x))=\deg(g_2(x))$. This implies that $a=0$, and hence that $\sigma$ is the identity, since otherwise $g_2(x)$ would contain a $\deg(h_2(x))$ term. Therefore the two standard forms must be equal.

We now prove (ii). By definition if $g(x),h(x)\in k[x]$, where $g(x)/h(x)$ is the standard form for $f(x)$, then $f(x)$ is equivalent to a rational function defined over $k$. For the other direction, suppose that $f(x)$ is equivalent to $f_k(x)$, a rational function defined over $k$. Applying the procedure from the first paragraph of this proof to produce the standard form of $f_k(x)$, we see that all fractional linear transformations are defined over $k$. This implies that the standard form for $f_k(x)$, which is also the standard form for $f(x)$, is defined over $k$.

Finally we show that having a standard form $g(x)/h(x)$ defined over $k$ is equivalent to the corresponding point in $G_d^{sep}$ having residue field $k$. Recall that $\gds$ is an open subscheme of $Gr(2,Poly_k(d))$, the Grassmannian of 2-planes in the space of polynomials of degree at most d. The equivalence class of $g(x)/h(x)$ is contained in the open subscheme of $Gr(2,Poly_k(d))$ consisting of equivalence classes of rational functions with a representative of the form $\tw{g}(x)/\tw{h}(x)$, where $\tw{g}(x)$ and $\tw{h}(x)$ are monic, $\deg(\tw{h}(x))<\deg(\tw{g}(x))$, and $\tw{g}(x)$ contains no $\deg(h(x))$ term. This subscheme is isomorphic to $\mathbb{A}^{2d-2}_k$. Under the natural isomorphism the unspecified coefficients of $\tw{g}(x)$ and $\tw{h}(x)$ serve as the coordinates of affine space, which establishes the desired equivalence.
\end{proof}

\section{Sufficient Conditions in Characteristic 3}
Roughly speaking, our strategy for proving \cref{mainthm} will be to show that there is no residue field extension for the map $Wr:\gds\rightarrow\PP^{2d-2}_k$ over points corresponding to polynomials with no repeated roots and with all coefficients contained in $k$. Our initial argument will imply that each such point has a unique preimage, but we must then rule out the case of purely inseparable field extensions by appealing to the following:

\begin{lem}\label{noinsep}
Let $k$ be a field of characteristic $p>2$ and $d$ a positive integer. The map $Wr:G_d^{sep}\rightarrow \mathbb{P}^{2d-2}_k$ is not generically purely inseparable (i.e.\ the induced extension of function fields is not purely inseparable).
\end{lem}

\begin{proof}
First note that by \cite[Theorem 1.4]{ossp} there is a non-zero number of equivalence classes of degree $d$ rational functions simply ramified at a generic choice of points $P_1,...,P_{2d-2}\in\PP^1_{\overline{k}}$. Since a generic polynomial has no repeated roots, this implies that $Wr$ is dominant, and hence induces an extension of function fields.

We shall explicitly describe the map $Wr$ in coordinates. Consider $c_1,...,c_{2d-2}$ as parameters on $\PP^1_{\overline{k}}$. Let $g(x)/h(x)$ be the standard form of a degree $d$ rational function where $g(x)=\sum a_i x^i$ and $h(x)=\sum b_i x^i$. If the Wronskian of $g(x)/h(x)$ maps to $\prod (x-c_i)$ (where if $c_i=\infty$ we define $x-c_i$ to equal 1), it must be the case that $h(x)g'(x)-g(x)h'(x)=\alpha\prod (x-c_i)$, where $\alpha$ is the leading coefficient of the left hand side. Expanding both sides and equating coefficients of $x^j$ for $0\leq j\leq 2d-1$, we have that $q_j(a_i,b_i)=(-1)^j e_{2d-2-j}(c_i)$, where $e_k(c_i)$ is the $k^{th}$ elementary symmetric polynomial in the $c_i$ and each $q_j(a_i,b_i)$ is a polynomial of total degree at most two.

The function field of the codomain of $Wr$ is the fixed field of $k(c_1,...,c_{2d-2})$ under the action of the symmetric group which permutes the $c_i$. It is well known that a $k$-basis for this fixed field is the symmetric polynomials $e_1(c_i),...,e_{2d-2}(c_i)$. Therefore the field extension induced by $Wr$ is given by $e_k\mapsto q_{2d-2-k}(a_i,b_i)$. Since the degree of each $q_j$ is at most two and there are no relations amongst the $a_i,b_i$, this field extension is not purely inseparable.
\end{proof}

Equipped with \cref{noinsep}, we now prove our main result:

\begin{thm}\label{mainthm}
Let $k$ be a field of characteristic 3. Every simply ramified rational function defined over $\overline{k}$ with only $k$-rational ramification points is equivalent to a rational function defined over $k$.
\end{thm}

\begin{proof}
Fix an integer $d>1$. By \cite[Theorem 1.4]{ossp}, there is exactly one equivalence class of degree $d$ rational functions simply ramified at a generic choice of points $P_1,...,P_{2d-2}\in\PP^1_{\overline{k}}$. Indeed, using the notation of \cite[Theorem 1.4]{ossp}, since the ramification is all simple we have that $e_1,...,e_{2d-2}=2$. The only possible $(2d-5)$-tuple of $e_i'$ satisfying conditions (i) and (ii) in \cite[Theorem 1.4]{ossp} begins with a 1 and alternates between 1 and 2, which implies that there is generically one such equivalence class, as claimed.

Switching back to the notation introduced in this paper, let $T$ be the $(2d-2)$-tuple consisting of all ones. Then $X_T$ is the space of equivalence classes of simply ramified degree $d$ rational functions. By \cite[Theorem 5.1]{myconj}, the number of equivalence classes of such rational functions is finite. In other words, $Wr:X_T\rightarrow\PP^{2d-2}_k$ is quasi-finite.

By Zariski's Main Theorem, we have a factorization

$$\xymatrix{
X_T \ar[rr]^{Wr}\ar[rd]_\iota &  & \PP_k^{2d-2} \\
& X\ar[ur]_\psi
}$$
where $\iota$ is an open immersion and $\psi$ is finite. Since $\psi$ is finite and generically the fibers of $Wr$ contain a single point, this implies that $\psi$ must either be degree 1 or purely inseparable. Since the latter possibility is eliminated by \cref{noinsep}, $\psi$ must be degree 1. Therefore there is no residue field extension over closed points for the map $Wr:X_T\rightarrow \PP^{2d-2}_k$. Hence, any point corresponding to an equivalence class of rational functions with only $k$-rational ramification points must have residue field $k$, which implies the result by \cref{form}(ii).
\end{proof}

\section{Insufficient Conditions}\label{negativesect}

We now turn our attention to considering conditions under which we cannot guarantee a rational function with only $k$-rational ramification points is equivalent to a rational function defined over $k$. Our first result establishes that \cref{mainthm} does not hold in characteristic greater than 3. \cref{char2cor} will establish an analogous result for characteristic 2.

\begin{prop}\label{finitefieldsimple}
Let $k$ be a finite field of characteristic $p>3$. There exists a simply ramified degree 3 rational function defined over $\overline{k}$ with only $k$-rational ramification points which is not equivalent to a rational function defined over $k$.
\end{prop}

\begin{proof}
Rather than repeating the arguments, we note that the results in \cite{fabershapiro} stated in characteristic 0 in fact hold for any field $k$ of characteristic greater than 3. In particular, using the notation from \cite{fabershapiro}, since the resultant of $x^3+ux^2$ and $(2u+3)x-(u+2)$ is $2(u+2)^2(u+1)^2$, the only values of $u$ for which $f_u(x)$ degenerates into a lower degree rational function remain $u=-2,-1$ when $\cha(k)>3$.

By \cite[Lemma 2.1]{fabershapiro} every degree 3 rational function defined over $\overline{k}$ ramified at $0,1,\infty$ is equivalent to a unique rational function of the form
$$f_u(x)=\frac{x^3+ux^2}{(2u+3)x-(u+2)},$$
where $u\in\overline{k}$. By factoring the Wronskian, we see that the last ramification point of $f_u(x)$ is $\varphi(u)=-(u^2+2u)/(2u+3)$.

By \cite[Proposition 2.3]{fabershapiro}, $f_u(x)$ is equivalent to a rational function defined over $k$ precisely when $u\in k$. Therefore to conclude the proof it suffices to show that there exists a $u\in\overline{k}\setminus k$ such that $\phi(u)\in k\setminus\{0,1\}$. Here we exclude 0 and 1 because we wish to produce a simply ramified rational function. Note that $\phi(-3),\phi(-1)=1$ and $\phi(0),\phi(-2)=0$. Since $k$ is a finite field, this implies that $\phi: k\rightarrow k$ cannot be surjective. Choose a $c\in k\setminus \phi(k)$. Solving $\phi(u)=c$ over $\overline{k}$ yields the desired rational function.
\end{proof}

Next we show that the situation is particularly bad when a ramification index is at least as large as the characteristic of the ground field (referred to as ``low characteristic'' in \cite{ossp}).

\begin{prop}\label{noshapiro}
Let $k$ be a non-algebraically closed field of characteristic $p>0$ and $f(x)$ a rational function defined over $\overline{k}$ with only $k$-rational ramification points. Suppose that $f(x)$ has a ramification point with ramification index at least $p$. Then there exist infinitely many equivalence classes of rational functions of the same degree and Wronskian as $f(x)$ which are not equivalent to a rational function defined over $k$.
\end{prop}

\begin{proof}
Let $d$ be the degree of $f(x)$. Since there exists a fractional linear transformation defined over $k$ sending $\infty$ to the point with ramification index at least $p$, we may without loss of generality assume that $\infty$ is this point. Let $g(x)/h(x)$ be the standard form for $f(x)$. Consider the rational functions
\begin{align*}
f_t(x)&=\frac{g(x)}{h(x)}+t\cdot x^p=\frac{g(x)+t\cdot x^p h(x)}{h(x)}
\end{align*}
with parameter $t\in \overline{k}$. By our assumptions, each $f_t$ is degree $d$ and its Wronskian is $h(x)g'(x)-g(x)h'(x)$, independent of $t$. By post-composing $f_t$ with an arbitrary fractional linear transformation, one verifies that $f_t$ and $f_{t'}$ are equivalent only if $t=t'$. This implies that these rational functions correspond to a one dimensional subscheme of $\gds$ on which $Wr$ is constant. Any point on this subscheme with residue field unequal to $k$ corresponds to an equivalence class of rational functions, none of which are defined over $k$ by \cref{form}(ii).
\end{proof}

When the characteristic of $k$ is small, \cref{noshapiro} gives us a great deal of information. As our first case, since all ramified points have ramification index at least 2, it informally tells us that the analog of the Shapiro Shapiro conjecture does not hold at all in characteristic 2:

\begin{cor}\label{char2cor}
Let $k$ be a non-algebraically closed field of characteristic 2. Suppose that $f(x)$ is a degree $d>1$ rational function defined over $\overline{k}$ with only $k$-rational ramification points. Then there are infinitely many equivalence classes of degree $d$ rational functions with the same Wronskian as $f(x)$ which are not equivalent to a rational function defined over $k$.
\end{cor}

Additionally, \cref{noshapiro} tells us that in characteristic 3 the analog of the Shapiro Shapiro conjecture does not hold when the ramification is not simple:
\begin{cor}\label{char3cor}
Let $k$ be a non-algebraically closed field of characteristic 3. Suppose that $f(x)$ is a degree $d$ rational function defined over $\overline{k}$ with only $k$-rational ramification points. If the ramification of $f(x)$ is not simple, then there are infinitely many equivalence classes of rational functions of the same degree and Wronskian as $f(x)$ which are not equivalent to a rational function defined over $k$.
\end{cor}
\bibliography{RyanEberhart}{}

\providecommand{\bysame}{\leavevmode\hbox to3em{\hrulefill}\thinspace}
\providecommand{\MR}{\relax\ifhmode\unskip\space\fi MR }
\providecommand{\MRhref}[2]{%
  \href{http://www.ams.org/mathscinet-getitem?mr=#1}{#2}
}
\providecommand{\href}[2]{#2}
\begin{thebibliography}{MTV09}

\bibitem[Ebe13]{myconj}
Ryan Eberhart, \emph{Families and moduli of covers of curves with fixed
  ramification}, {\tt arXiv:1312.7144 [math.AG]} (2013).

\bibitem[EG02]{shapiro}
A.~Eremenko and A.~Gabrielov, \emph{Rational functions with real critical
  points and the {B}. and {M}. {S}hapiro conjecture in real enumerative
  geometry}, Ann. of Math. (2) \textbf{155} (2002), no.~1, 105--129.
  \MR{1888795}

\bibitem[FT16]{fabershapiro}
Xander Faber and Bianca Thompson, \emph{A very elementary proof of the {B.} and
  {M.} {S}hapiro conjecture for cubic rational functions},
  arXiv:math.NT/1601.04938 (2016).

\bibitem[MTV09]{fullshapiro}
Evgeny Mukhin, Vitaly Tarasov, and Alexander Varchenko, \emph{The {B}. and {M}.
  {S}hapiro conjecture in real algebraic geometry and the {B}ethe ansatz}, Ann.
  of Math. \textbf{170} (2009), no.~2, 863--881.

\bibitem[Oss06]{ossp}
Brian Osserman, \emph{Rational functions with given ramification in
  characteristic $p$}, Compositio Math. \textbf{142} (2006), no.~2, 433--450.

\end{thebibliography}
\bibliographystyle{amsalpha}
\end{document}